\documentclass[a4paper]{amsart}
\usepackage{amsfonts}
\usepackage{amssymb,xspace}
\usepackage{xy}
\let\tilde\widetilde
\let\hat\widehat
\let\wt\widetilde
\let\ov\overline
\let\<\langle
\let\>\rangle
\DeclareMathOperator{\ad}{ad} \DeclareMathOperator{\GL}{GL}
\DeclareMathOperator{\Id}{Id} \DeclareMathOperator{\diag}{diag}
\def\cA{\mathcal{A}}
\def\cE{\mathcal{E}}
\def\cH{\mathcal{H}}
\def\cO{\mathcal{O}}
\def\cQ{\mathcal{Q}}
\def\cT{\mathcal{T}}
\def\cU{\mathcal{U}}
\def\cV{\mathcal{V}}
\def\CC{\mathbb{C}}
\def\ZZ{\mathbb{Z}}
\renewcommand{\MR}[1]{}

\numberwithin{equation}{section}

\theoremstyle{plain}
\newtheorem{theorem}[equation]{Theorem}
\newtheorem{proposition}[equation]{Proposition}
\newtheorem{lemma}[equation]{Lemma}
\newtheorem{corollary}[equation]{Corollary}
\newtheorem*{claim*}{Claim}
\newtheorem{claim}{Claim}

\theoremstyle{definition}

\theoremstyle{remark}
\newtheorem{remark}[equation]{Remark}

\begin{document}

\title[Frobenius manifolds with tt*-structure]
{Flat meromorphic connections of Frobenius manifolds with
tt*-structure}

\author[J.~Lin]{Jiezhu Lin}
\address[J.~Lin]{Guangzhou, China}
\email{ljzsailing@163.com}

\author[C.~Sabbah]{Claude Sabbah}
\address[C.~Sabbah]{UMR 7640 du CNRS\\
Centre de Math\'ematiques Laurent Schwartz\\
\'Ecole polytechnique\\
F--91128 Palaiseau cedex\\
France} \email{sabbah@math.polytechnique.fr}
\urladdr{http://www.math.polytechnique.fr/~sabbah}

\thanks{C.S.\ was supported by the grant ANR-08-BLAN-0317-01 of the Agence nationale de la recherche.}

\begin{abstract}
The base space of a semi-universal unfolding of a hypersurface
singularity carries a rich geometric structure, which was
axiomatized as a CDV-structure by C.~Hertling. For any
CDV-structure on a Frobenius manifold~$M$, the pull-back bundle
$\pi^* \mathcal{T}_M^{(1, 0)}$ by the projection
$\pi:\mathbb{C}\times M\to M$ carries two natural holomorphic
structures equipped with two flat meromorphic connections. We show
that, for any semi-simple CDV-structure, there is a formal
isomorphism between these two bundles compatible with connections.
Moreover, if we assume that the super-symmetric index $\mathcal Q$
vanishes, we give a necessary and sufficient condition for such a
formal isomorphism to be convergent, and we make it explicit for
$\dim M = 2$.
\end{abstract}


 \maketitle

\setcounter{section}{-1} \vspace*{-\baselineskip}

\section{Introduction}\label{section0}
tt* geometry, which appeared first in papers of Cecotti and Vafa
(\cite{CV3}, \cite{CVN}), is now understood after the work of C.
Hertling \cite{Hert2}, as an enrichment of that of harmonic Higgs
bundle $(E, \Phi, h)$ previously introduced by N. Hitchin and C.
Simpson. The one-to-one correspondence between harmonic Higgs
bundles and variations of polarized twistor structures of weight
$0$ can be extended, after the work of C.~Hertling \cite{Hert2}
(respectively, C.~Sabbah \cite[Chap.~7]{Sabb08}) to that between
CV-structures (respectively, integrable harmonic Higgs bundles)
and variations of pure TERP structures, cf.~\cite{Hert2} for the
terminology (respectively, integrable variations of Hermitian pure
twistor structures of weight $0$, cf.~\cite{Sabb08} for the
terminology). The important object of a variation of pure TERP
structure (respectively, integrable variations of Hermitian pure
twistor structures of weight $0$) is the twistor bundle on
$\mathbb{P}^1 \times M$ with a flat $C^\infty$-connection whose
holomorphic part has a pole with Poincar\'e rank one along $\{ 0
\} \times M$ and whose anti-holomorphic part has a pole with
Poincar\'e rank one along $\{ \infty \} \times M$. Such a
connection will be called the structure connection of the
CV-structure.

Of particular interest for us is the case where a CV-structure
exists on a Frobenius manifold $(M, g, \circ, e, \mathcal{E})$ in
a compatible way. Such a structure is called a CDV-structure, and
is defined and studied by C. Hertling in \cite{Hert2}. In such a
case, there are two natural holomorphic structures on the
pull-back tangent bundle $\pi^* \mathcal{T}_M^{(1,0)}$ by
$\pi:\mathbb{C} \times M\to M$, each of which carries a flat
connection. More precisely, one is the holomorphic vector bundle
$\mathcal{H}_1:=\pi^* \mathcal{T}_M$ with the structure connection
of the Frobenius manifold $M$, denoted by $\widetilde{\nabla}$,
which is an integrable meromorphic connection and has Poincar\'e
rank one along $\{ 0 \} \times M$ (and is extended as a
logarithmic connection along $\{ \infty \} \times M$); the other
one is the structure connection of the CV structure, denoted by
$\widetilde{D}$, on $\pi^* \mathcal{T}_M^{(1,0)}$. This connection
is integrable, hence the $(0,1)$-part of this connection, which
has no pole on $\mathbb{C} \times M$, gives another holomorphic
structure, denoted by $\mathcal{H}_2$, on $\pi^*
\mathcal{T}_M^{(1,0)}$. Moreover, $\widetilde{D}$ is a meromorphic
connection on $\mathcal{H}_2$ having Poincar\'e rank one along $\{
0 \} \times M$. We will say that the CDV-structure is
\emph{strongly potential} if both holomorphic bundles with
connection $(\cH_1,\wt\nabla)$ and $(\cH_2,\wt D)$ are
\emph{isomorphic} (and an isomorphism between both is called a
potential). A CDV-structure is strongly potential if and only if
there exists a holomorphic isomorphism
\begin{equation*}
\phi: \mathcal{H}_1 \longrightarrow \mathcal{H}_2,
\end{equation*}
such that
\begin{equation}\label{eq:Cphi}
\phi \widetilde{\nabla} = \widetilde{D} \phi.
\end{equation}

In \cite[Th.~5.15]{Hert2}, C.~Hertling gives a criterion to
produce CDV-structures, and the resulting structures are strongly
potential. The terminology used here comes from \cite{Sabb08}
where the notion of a potential harmonic Frobenius manifold is
considered, and where the criterion of Hertling is shown to
produce a potential harmonic Frobenius manifold with stronger
properties. It should be noticed that the CDV-structure
constructed by Hertling on the base space of the universal
unfolding of a hyper-surface singularity (cf.~\cite[\S8]{Hert2}),
and that constructed by Sabbah for the universal unfolding of a
convenient and non-degenerate Laurent polynomial
(cf.~\cite[\S4.c]{Sabb08}) both use Hertling's criterion, and
therefore give rise to a strongly potential CDV-structure
(however, in \cite[\S8]{Hert2} this is not mentioned).

The purposes of this article are to construct a formal isomorphism
between these two bundles with connections for any semi-simple
CDV-structures, and to analyze the strength of the potentiality
property in simple examples of CDV-structures. As~a~part of the
data of a CDV-structure is a self-adjoint operator $\cQ$ on
$\mathcal{T}_M^{(1,0)}$ (the ``new super-symmetric index''). More
generally, such an operator exists on the underlying bundle $K$ of
a CV-structure (cf.~\cite{Hert2} and \cite[Chap.\,7]{Sabb08}). If
the CV-structure corresponds to a variation of polarized Hodge
structures, then the eigenvalue decomposition of~$\cQ$ corresponds
to the Hodge decomposition, and the eigenvalues correspond to the
Hodge exponents. The simplest examples of variation of polarized
Hodge structures are those of Tate type, of pure type $(0,0)$.
They are nothing but flat Hermitian vector bundles. By analogy, we
will call a \emph{Tate CV-structure} (resp.\ a \emph{Tate
CDV-structure}) a CV-structure (resp.~a~CDV-structure) such that
$\cQ=0$. In this article, we will restrict our attention to Tate
CDV-structures when analyzing the strength of the potentiality
property. The existence of a Tate CDV-structure on any semi-simple
Frobenius manifold has been discussed in \cite{Lin09}, where the
author gave explicitly the Hermitian metric and proved that such a
CDV-structure is a harmonic Frobenius manifold (cf.~\cite{Sabb08}).

For many interesting examples of CDV-structures, the two
meromorphic connections $(\cH_1,\wt\nabla)$ and $(\cH_2,\wt D)$
have irregular singularities along $\{0\}\times M$, hence one
cannot hope in general that there exists a holomorphic isomorphism
between these two bundles compatible with connections.

The main result (Theorem \ref{Q0m}) of this article is to show
that, for any semi-simple CDV-structure on a $1$-connected complex
manifold $M$, a \emph{formal} isomorphism between the holomorphic
bundles $\mathcal{H}_1$ and $\mathcal{H}_2$ compatible with the
meromorphic connections does exist, and we notice that a
holomorphic lift of this isomorphism exists as soon as it exists
for the restriction at one point of $M$ of the bundles with
connection. We are then able to give a necessary and sufficient
condition for such a formal isomorphism to be convergent for a
Tate CDV-structure, and in the case that $\dim M=\nobreak2$ we
give an explicit formula (Corollary \ref{Q02}) for such an
isomorphism.

The proof of Theorem \ref{Q0m} consists in constructing a formal
isomorphism $\hat{\phi}{}^{o}$ between two restricted bundles
$\mathcal{H}_i |_{\mathbb{C} \times o}$ compatible with the
corresponding restricted connections for every point $o \in M$.
Since we assume that the underlying Frobenius manifold is
semi-simple, we can apply a theorem of Malgrange \cite{Malg}
(cf.~also \cite[Th. II.2.10]{Sabb}) to extend it and get a formal
isomorphism $\hat{\phi}$ between the bundles $\mathcal{H}_i$
compatible the meromorphic connections. According to the local
constancy of the Stokes sheaf, we also obtain that $\hat\phi$
lifts to a holomorphic isomorphism compatible with the connections
if and only if $\hat\phi{}^o$ does so.

\subsubsection*{Notation and terminology}
We usually refer to \cite{D,Hert, Mani,Sabb} for the notion of
Frobenius manifold, that we denote by $(M, \circ, g, e,
\mathcal{E})$, where~$M$ is a complex manifold, $g$~is a metric on
$M$ (that is, a symmetric, non-degenerate bilinear form, also
denoted by $\< \,, \>$) with associated Levi-Civita connection
denoted by $\nabla$, $\circ$ is a commutative and associative
product on $TM$ which depends smoothly on $M$, $e$ is the unit
vector field for $\circ$ and $\cE$ is called the Euler vector
field. These data are subject to conditions that we do not repeat
here.

In the following, we will restrict to semi-simple Frobenius
manifolds, i.e., the multiplication by $\cE$ has pairwise distinct
eigenvalues at each point of $M$, and we will assume the existence
of canonical coordinates, that we will denote by
$\boldsymbol{u}=(u^1,\dots,u^m)$, such that $\mathcal{E}= \sum_{k}
u^k \partial_{u^k}$. We note that, by assumption, we have $u^i
\neq u^j$ on $M$ if $i \neq j.$ We will set $e_i:=\partial_{u^i}$.
Such coordinates exist (at least in the \'etale sense) whenever
$M$ is $1$-connected (cf.~\cite[\S VII.1.8]{Sabb}), and we will
mainly restrict to this case. We then define $\cV_{ij}$ by
\[
\cV e_i = \Big(\nabla \mathcal{E}+ \frac{2-d}{2}\Id\Big)e_i :=
\sum_j \cV_{ij} e_j, \quad \forall i.
\]
Recall also that there exists a metric potential $\eta$, that is,
a function such that $\eta_i:=\partial\eta/\partial
u^i=g(e_i,e_i)$ for all $i$. The matrix $(\cV_{ij})$ can be
expressed in terms of this potential as
(cf.~\cite[Eq.~(87)]{Lin09})
\begin{equation}\label{eq:V}
\cV_{ij}=
\begin{cases}0&\text{if $i=j$},\\
(u^j-u^i)\partial_{u^i}\partial_{u^j}(\eta)/2\partial_{u
^j}(\eta)&\text{if $i\neq j$}.
\end{cases}
\end{equation}

For the notion and notation relative to TERP structures and
CDV-structures, we refer to \cite{Hert2}. For the deformations of
connections with poles of Poincar\'e rank one and their formal
decompositions, we refer to \cite[\S III.2]{Sabb} In particular,
we will use the real structure~$\kappa$ on $\cT_M^{(1,0)}$, which
determines, together with $g$, a Hermitian form $h$ on
$\cT_M^{(1,0)}$ whose associated Chern connection is denoted by
$D=D'+\overline\partial$. The self-adjoint operator $\cQ$ on
$\mathcal{T}_M^{(1,0)}$ is defined by $\mathcal{Q}:=
D_{\mathcal{E}}- \mathcal{L}_{\mathcal{E}}- \frac{2-d}{2}\cdot
\Id.$ If the given CDV-structure is semi-simple, we define
$\mathcal{Q}_{ij}$ by
\[
\cQ e_i = \Big( D_{\mathcal{E}}- \mathcal{L}_{\mathcal{E}}-
\frac{2-d}{2}\cdot \Id \Big)e_i =: \sum_j \cQ_{ij} e_j, \quad
\forall i.
\]
By straightforward computation, we get that
\begin{equation}\label{eq:Q}
\cQ_{ij}=
\begin{cases}0&\text{if $i=j$},\\
\omega_i^j (\cE)&\text{if $i\neq j$}.
\end{cases}
\end{equation}
Here $\omega_i^j$ are the connection forms of $D'$ under the local
frame $e_i$.

On the product $\CC\times M$, we usually denote by $z$ the
coordinate on $\CC$.

\subsubsection*{Acknowledgements}
The authors thank Claus Hertling for attracting their interest to
the question considered in this article, and for useful
discussions and comments.

\section{Main results}\label{section2}
Before we state the main theorem, we give some equivalent
conditions for a semi-simple Tate CDV$\oplus$-structure.

\begin{proposition}\label{propA}
Let $(M, g, \circ, e, \mathcal{E}, \kappa)$ be a semi-simple
CDV$\oplus$-structure with canonical coordinates as above. Let
$\eta$ be the associated metric potential. Then the following
statements are equivalent:
\begin{enumerate}
\item[a)] $(h(e_i, e_j))_{m \times m} = \diag(|\eta_1 |,\dots,
|\eta_m|)$;

\item[b)] $D'$ is a holomorphic connection, i.e.
$D'\overline{\partial}+ \overline{\partial} D'=0$;

\item[c)] $\mathcal{Q}=0$.
\end{enumerate}
\end{proposition}

\begin{theorem}\label{Q0m}
For any semi-simple CDV-structure on a $1$-connected complex
manifold $M$, there exists a formal isomorphism $\hat{\phi}$
between the formalized bundles with connections $(\hat{\cH_1},
\hat{\widetilde{\nabla}})$ and $(\hat{\cH_2},
\hat{\widetilde{D}})$, where we set
$\hat\cH=\hat\cO_M\otimes_{\cO_{\CC\times M}}\cH$ and
$\hat\cO_M=\varprojlim_k\cO_{\CC\times M}/z^k$. Moreover,
$\hat\phi$ lifts as a holomorphic isomorphism if and only if its
restriction $\hat\phi{}^o$ at one point $o\in M$ induces a
holomorphic isomorphism between the restricted holomorphic bundles
with connection. Lastly, if we assume that the CDV-structure is
Tate and positive, then such a holomorphic lift exists if and only
if the monodromy of $(\cH_1^o,\wt\nabla{}^o)$ is equal to
identity.
\end{theorem}

The following theorem will explain what happens when the
semi-simple Tate CDV$\oplus$-structure is strongly potential.

\begin{theorem}\label{thmA}
A semi-simple Tate CDV$\oplus$-structure on a $1$-connected
complex manifold $M$ with canonical coordinates as above is
strongly potential if and only if, for some point $o\in M$ with
canonical coordinates $(u^1_o,\dots,u^m_o)$, there exists a matrix
$(\psi_{ij}^o)\in \GL_{m \times m}(\CC[z])$ such that
\begin{equation}\tag*{(\protect\ref{thmA})$(*)$}\label{eq:thmA}
\partial_z (\psi_{ij}^o)=\dfrac{1}{z^2} (u^i_o-u^j_o)\psi_{ij}^o - \dfrac{1}{z}
\sum_{k}\cV_{ik}^o \psi_{kj}^o, \quad \forall i,j.
\end{equation}
The matrix $(\phi_{ij}^o)$ of an isomorphism $\phi^o$ satisfying
\eqref{eq:Cphi} at $o$ is then given by
\[
(\phi_{ij}^o)=(\psi_{ij}^o)\cdot\diag\big(\exp(z\ov
u{}^1_o),\dots, \exp(z\ov u{}^m_o)\big).
\]
\end{theorem}

\begin{corollary}\label{CoFlka}
Any semi-simple CDV$\oplus$-structure on a $1$-connected complex
analytic manifold $M$, such that $D'= \nabla$ (or equivalently,
such that the canonical coordinates are $\nabla$-flat), is of Tate
type and strongly potential. Moreover, in the frame $e_i$, the
isomorphism $\phi: (\mathcal{H}_1, \widetilde{\nabla}) \rightarrow
(\mathcal{H}_2, \widetilde{D})$ is determined by the following
matrix
\begin{equation*}
(\phi_{ij}(z, \boldsymbol{u}))_{m\times m} =
\diag(c_1\cdot\exp(z\ov u{}^1),\dots, c_m \cdot\exp(z\ov u{}^m)),
\end{equation*}
for some nonzero constants $c_1,\dots, c_m$.
\end{corollary}

When $\dim M =2$, we will make explicit the necessary and
sufficient condition for a semi-simple Tate CDV$\oplus$-structure
on $M$ to be strongly potential, as given in Theorem \ref{Q0m}.

\begin{corollary}\label{Q02}
Let $M$ be a $1$-connected complex analytic manifold $M$ with
$\dim M =2$. Let $(M, \circ, e, \mathcal{E}, \kappa)$ be a
semi-simple Tate CDV$\oplus$-structure on $M$ such that
$g(e,e)\neq0$. Let~$d$ be the constant such that $\nabla
\mathcal{E} + (\nabla \mathcal{E})^* =(2-d)\cdot \Id$. Then the
CDV$\oplus$-structure is strongly potential if and only if
$d\in2\ZZ$.
\end{corollary}

\section{Proof of the theorems}\label{section3}

We will use a system of holomorphic canonical coordinates $u^1,
u^2,\dots,u^m$ of the Frobenius manifold as above and we set $e_i
= \partial_{u^i}$. We will then use the following notation:
\begin{align}
\kappa (e_i)&= \sum_{k} K_{ik}e_k,\label{eq:Mkappa}\\
e_i\circ e_j=-\Phi_{e_i} (e_j)&= \sum_{k} {C^{(i)}}_{j}^{k}e_k,\label{eq:MPhi}\\
-\Phi^{\dag}_{ \overline e_i} (e_j)&=
\sum_{k}\tilde{C^{(i)}}_{j}^{k} e_k, \quad\text{$\Phi^\dag$ is the
$h$-adjoint of $\Phi$}. \label{eq:MPhidag}
\end{align}
By definition of the canonical coordinates, we have
${C^{(i)}}_{j}^{k}= \delta_{ik} \cdot \delta_{jk}$. Because of
$h(X, Y)= g(X, \kappa Y)$ and $\Phi^{*}=\Phi$, we have
$\Phi^{\dag}= \kappa \Phi \kappa$ since, for all $X,Y$,
\[
h(X, \Phi^{\dag}Y)=h(\Phi X, Y)= g(\Phi X, \kappa Y)= h(X, \kappa
\Phi \kappa Y).
\]
This is expressed by $\tilde{C^{(i)}}= \overline{K} \cdot\nobreak
\overline{C^{(i)}} \cdot\nobreak K$.

\begin{proof}[Proof of proposition \ref{propA}]\mbox{}
\par\nopagebreak\noindent
a) $\Rightarrow$ b) and a) $\Rightarrow$ c). This is proved in
\cite[Th.\,2]{Lin09}.

\noindent b) $\Rightarrow$ a). If $D'$ is holomorphic, then by the
harmonicity condition (cf.~\cite[(1.1)]{Sabb08}), we have
\begin{equation*}
\Phi \wedge \Phi^\dag + \Phi^\dag \wedge \Phi =0.
\end{equation*}
By straightforward computations we get
\begin{equation}\label{eq:CC}
[C^{(i)}, \overline{K} \overline{C^{(j)}} K]=0, \quad\forall i, j.
\end{equation}
Computing \eqref{eq:CC} directly, we conclude that, for any $i$
there exists a unique $j_i$ such that
$$K_{i j_i} \neq 0,\quad K_{j_ii} \neq 0.$$
However, $$h_{ii}= K_{ii} \cdot \eta_{i} > 0.$$ Hence we get that
\begin{equation}\label{eq:dK}
K=\diag(K_{11},\dots, K_{mm}).
\end{equation}
The relations $\kappa^2 = \Id$ and \eqref{eq:dK} imply
$$|K_{ii}| = 1, \quad\forall i.$$
So we can conclude that
\begin{equation*}
h_{ii}= |\eta_{i}|, \quad\forall i.
\end{equation*}

c) $\Rightarrow$ b). Assume that $\cQ=0$. From Equations
(1.12)-(1.16) in \cite{Sabb08} together with the $h$-adjoint ones,
we obtain
\begin{equation}\label{eq:UC0}
D(\cU)=D'(\cU) =- \Phi, \quad [\mathcal{U}, D(\mathcal{U})]=0,
\quad [\mathcal{U}, D(\mathcal{U}^\dag)]=0
\end{equation}
and adjoint equations. Then $^V\!D:= D+\Phi + \Phi^\dag$ can be
written as $D- D(\mathcal{U}+ \mathcal{U}^\dag)$. Since
$D(\mathcal{U}+ \mathcal{U}^\dag)$ commutes with $\mathcal{U}$ and
$\mathcal{U}^\dag$ by \eqref{eq:UC0}, it commutes with
$\mathcal{U}+\mathcal{U}^\dag$ and we have
$$^V\!D = e^{\mathcal{U}+\mathcal{U}^\dag} D e^{-(\mathcal{U}+\mathcal{U}^\dag)}.$$
Recall now that, for a CDV$\oplus$ structure, the connection
$^V\!D$ is flat, being the restriction to $z=1$ of the flat
connection $\wt D+(\ov\partial+z\Phi^\dag)$. This implies
therefore that~$D$ is flat, in particular $D'$ is a holomorphic
connection.
\end{proof}

\begin{remark}\label{rem:phidag}
If the properties of the proposition are satisfied, then
\[
\wt{C^{(i)}}=\ov{C^{(i)}}={C^{(i)}},\quad\forall i.
\]
\end{remark}

\begin{corollary}\label{dUdag}
Let $(M, g, \circ, e, \mathcal{E}, \kappa)$ be a semi-simple Tate
CDV$\oplus$-structure with canonical coordinates as above and let
$\eta$ be the associated metric potential. Then
\begin{enumerate}
\item[1)] $(\mathcal{U}^{\dag}_{ij})=
\overline{\mathcal{U}}=\diag(\ov u{}^1,\dots, \ov u{}^m)$, where
$(\mathcal{U}^{\dag}_{ij})$ is the matrix of $\mathcal{U}^{\dag}$
in the local frame $e_i:= \partial_{u^i}$. In particular,
$[\mathcal{U}, \mathcal{U}^\dag]=0$;

\item[2)] $D'=\nabla$ if and only if the canonical local
coordinates $u^i$ are $\nabla$-flat, i.e.,
$$
\nabla \partial_{u^i}=0, \quad \forall i.\eqno\qed
$$
\end{enumerate}
\end{corollary}

\begin{proof}[Proof of Theorem \ref{Q0m}]
Given any semi-simple CDV-structures on $M$, given a point $o \in
M$, let $(H_1, \hat{\nabla})$ (resp.~$(H_2, \hat{D})$) denote by
the restricted bundles with connections at $o$ of $(\hat{\cH_1},
\hat{\widetilde{\nabla}})$ (resp.~$(\hat{\cH_2},
\hat{\widetilde{D}})$). We will firstly prove that there exists a
formal isomorphism between $(H_1, \hat{\nabla})$ and $(H_2,
\hat{D})$. For simplicity, we denote by~$U$ (resp.~$V$ and $Q$)
the restriction at $o$ of $\mathcal{U}$ (resp.~$\mathcal {V}$,
$\mathcal {Q}$). We note that $V$ (resp.~$Q$) is in the image of
$\ad U$: this can be seen by considering the local frame $s_i:
=(\pi^*e_i)|_{\mathbb{C}\times o}$ of $H_1$ (resp.~$H_2$),
according to the relation \eqref{eq:V} (resp.~\eqref{eq:Q}). Then
the system $(U- zV) dz/z^2$ (resp.~$(U - zQ - z^2 U^\dag) dz/z^2$)
is equivalent, by a holomorphic base change, to a system $(U + z^2
C(z)) dz/z^2$ (resp.~$(U + z^2 \widetilde{C}(z)) dz/z^2$), where
$C(z)$ (resp.~$\widetilde{C}(z)$) is a matrix whose entries are
holomorphic functions. The following lemma will imply that $(U +
z^2 C(z)) dz/z^2$ and $(U + z^2 \widetilde{C}(z)) dz/z^2$ are
formally isomorphic to $Udz/z^2$, therefore we conclude that they
are formally isomorphic.

\begin{lemma}\label{Qneq0}
Given any matrix $\hat C(z)$ whose entries are formal series
w.r.t.~$z$, then the system with matrix $(U + z^2 \hat C(z))
dz/z^2$ is equivalent, by formal base change, to a system
$Udz/z^2$ if we assume that $U$ is a regular semi-simple matrix.
\end{lemma}

\begin{proof}[Proof of lemma \ref{Qneq0}]
According to \cite[Th.\,III.2.15]{Sabb}, since $U$ is regular
semi-simple, the connection matrix $(U + z^2 \hat C(z)) dz/z^2$ is
equivalent to a unique diagonal matrix, each diagonal term being
equal to $(u_i+\mu_iz)dz/z^2$, where $u_i$ are the eigenvalues of
$U$ and $\mu_i\in\CC$. However, the coefficient of $dz/z$ in $(U+
z^2 \hat C(z)) dz/z^2$ is zero, hence by straightforward
computation, we conclude that all the constants $\mu_i$ are zero.
Therefore, $(U+ z^2 \hat C(z)) dz/z^2$ is formally equivalent to
$Udz/z^2$.
\end{proof}

Let us continue our proof. We have proved that for any point $o\in
M$ the restricted bundles with connections are formally
isomorphic. The existence of a formal isomorphism between the
bundles themselves is given by the following lemma.

\begin{lemma}\label{extend}
Let $(E, \nabla)$ and $(F, D)$ be two holomorphic bundles on
$\mathbb{C}\times M$ with meromorphic connections of Poincar\'e
rank one along $\{0 \} \times M$ whose ``residues'' $R^\nabla_0
(x)$ and $R^D_0(x)$ are regular semi-simple and have the same
eigenvalues $\lambda_i(x)$, where $M$ is a $1$-connected complex
manifold. Let $(\hat{E}, \hat{\nabla})$ and $(\hat{F}, \hat{D})$
denote by the formal bundles with connections associated to $(E,
\nabla)$ and $(F, D)$. If there exists a formal isomorphism
between the restricted bundles with connections of $(\hat{E},
\hat{\nabla})$ and $(\hat{F}, \hat{D})$ at one point $o \in M$,
then we can extend it as a formal isomorphism between $(\hat{E},
\hat{\nabla})$ and $(\hat{F}, \hat{D})$.
\end{lemma}

\begin{proof}[Proof of lemma \ref{extend}]
Since $R^\nabla_0 (x)$ (resp.~$R^D_0(x)$) is regular semi-simple
and has eigenvalues $\lambda_i(x)$, we know by
\cite[Th.\,III.2.15]{Sabb}  that $(\hat{E}, \hat{\nabla})$
(resp.~$(\hat{F}, \hat{D})$) can be decomposed in a unique way as
a direct sum of subbundles with connections of rank one. Moreover, if we assume
that $M$ is $1$-connected, then $\hat{E}$ (resp.~$\hat{F}$) is
trivializable and admits a basis in which the matrix of connection
forms of $\hat{\nabla}$ (resp.~$\hat{D}$) can be written as
$\diag(\omega^\nabla_1, \dots, \omega^\nabla_m)$ (resp.
$\diag(\omega^D_1, \dots, \omega^D_m)$), where $\omega^\nabla_i$
(resp.~$\omega^D_i$) takes the form $-d(\lambda_i(x) /z) +
\mu^\nabla_i dz/z$ (resp.~$-d(\lambda_i(x) /z) + \mu^D_i dz/z$),
and $\mu^\nabla_i$ (resp.~$\mu^D_i$) are some complex number. If
there exists a formal isomorphism between the restricted bundles
with connections of $(\hat{E}, \hat{\nabla})$ and $(\hat{F},
\hat{D})$ at one point $o \in M$, then the uniqueness of the
decomposition implies that $\mu^\nabla_i = \mu^D_i$ for every $i$.
Therefore $(\hat{E}, \hat{\nabla})$ and $(\hat{F}, \hat{D})$ are
formally isomorphic.
\end{proof}

Assume now that the restriction at a point $o\in M$ with
coordinates $(u^1_o,\dots,u^m_o)$ of the Saito connection is
holomorphically isomorphic to the connection with matrix $(U/z - Q
-z U^\dag)dz/z$. Recall (cf.~\cite[\S II.6]{Sabb}) that the
meromorphic Saito connection $(\cH_1[1/z],\wt\nabla)$
(resp.~$(\cH_2[1/z],\wt D)$) can be reconstructed, up to
isomorphism, from its formalization along $z=0$ together with a
section $\sigma_1$ (resp.~$\sigma_2$) of the Stokes sheaf. Since
$M$ is simply connected and since the Stokes sheaf is locally
constant (cf.~\cite[Th.\,II.6.1]{Sabb}), such a section is
constant and uniquely determined by its germ at $o$. According to
the base change property \cite[Prop.\,II.6.9]{Sabb}, this germ at
$o$ is the Stokes cocycle of the restricted bundle with
connection. Our assumption is that the germs at $o$ of $\sigma_1$
and $\sigma_2$ are equal. Therefore, the section $\sigma_1$ is
equal to $\sigma_2$, that is, the formal isomorphism between the
meromorphic bundles $(\cH_1[1/z],\wt\nabla)$ and $(\cH_2[1/z],\wt
D)$ is convergent. Since it induces a formal isomorphism between
the holomorphic bundles $(\cH_1,\wt\nabla)$ and $(\cH_2,\wt D)$,
it is indeed a holomorphic isomorphism between these bundles.

Let us prove the last statement in the theorem. Given any Tate
semi-simple CDV$\oplus$-structure, we consider the bundle
$\pi^*\cT_M^{(1,0)}$ with the holomorphic structure
$\overline{\partial}+ z \Phi^\dag$, equipped with the meromorphic
connection
\begin{equation*}
\widetilde{D}= D' +d'_z+ \frac{1}{z} \Phi + \Big(\frac{1}{z}
\mathcal{U}-z\mathcal{U}^\dag\Big) \frac{dz}{z}.
\end{equation*}
We can write
\begin{equation}\label{eq:Phidag}
\overline{\partial}+ z \Phi^\dag = \overline{\partial}- z
\overline{\partial}(\mathcal{U}^\dag) = e^{z \mathcal{U}^\dag}
\circ \overline{\partial} \circ e^{-z \mathcal{U}^\dag}.
\end{equation}
Therefore, this holomorphic bundle $\mathcal{H}_2$ with
$z$-meromorphic connection $\widetilde{D}$ is isomorphic to the
holomorphic bundle $\mathcal{H}_1:=\ker\overline\partial$ equipped
with
$$e^{-z \mathcal{U}^\dag} \circ \Big[D' +d'_z+ \frac{1}{z} \Phi + \Big(\frac{1}{z} \mathcal{U}-z\mathcal{U}^\dag\Big) \frac{dz}{z}\Big] \circ e^{z \mathcal{U}^\dag},$$
which is written as
\begin{equation}\label{eq:ezu}
D' +d'_z+ \frac{1}{z} \Phi + e^{-z \mathcal{U}^\dag} \mathcal{U}
e^{z \mathcal{U}^\dag} \cdot \frac{dz}{z^2}.
\end{equation}
However, by Corollary \ref{dUdag}, we know that $\mathcal{U}^\dag$
commutes with $\mathcal{U}$, hence the connection in
\eqref{eq:ezu} can be written as
\begin{equation}\label{eq:wD}
D' +d'_z+ \frac{1}{z} \Phi + \mathcal{U} \cdot \frac{dz}{z^2}
\end{equation}
Let us begin with the connection given by \eqref{eq:wD}. Choose
canonical local coordinates $u^1,\dots, u^m$ of the underlying
semi-simple Frobenius manifold. Set $S_i := \pi^* \partial_{u^i}$,
for all~$i$. Then $(S_i)_{i=1,\dots,m}$ is a holomorphic local
frame for $\mathcal{H}_1$. Obviously, $\Phi$ is diagonal in this
frame, and $S_j$ is the eigenvector of $\mathcal{U}$ with
eigenvalue $u^j$. It follows that the connection \eqref{eq:wD} is
holomorphically isomorphic to the direct sum of the connections
$d'+ d'_z - (d'+ d_z')(u^j/z)$ for all $j=1, 2,\dots,m$. In
particular, for each such bundle, the monodromy around $z=0$ is
equal to the identity. Hence the existence of a holomorphic lift
will imply that the monodromy of $(U -z V)dz/z^2$ is equal to
identity.

On the other hand, the structure connection $\widetilde{\nabla}$
on the holomorphic bundle $\mathcal{H}_1$, under the above
holomorphic frame $S_i$, can be reduced to
\begin{equation}\label{eq:wnabla}
\nabla +d'_z+ \frac{1}{z} \Phi + (\mathcal{U}- z \cV) \cdot
\frac{dz}{z^2},
\end{equation}
where $\cV$ is given by \eqref{eq:V}.

We will use the following lemma to prove the other direction of
the last statement.
\begin{lemma}\label{lem:UV}
Consider a system $(U-zV)dz/z^2$, where $U$ is diagonal with
pairwise distinct eigenvalues. The following properties are
equivalent:
\begin{enumerate}
\item\label{enum:1} the system is meromorphically (resp.\
holomorphically) isomorphic to the system with matrix $Udz/z^2$,
\item\label{enum:2} the monodromy is equal to the identity (resp.\
and the diagonal part $V_{\diag}$ is zero).
\end{enumerate}
If these properties are satisfied, $V$ is semi-simple with
integral eigenvalues and integral (resp.\ zero) diagonal part.
\end{lemma}
To conclude the proof of Theorem \ref{Q0m}, it is enough,
according to \eqref{eq:wnabla}, to apply the lemma with $U=\cU^o$,
$V=\cV^o$, since we know by \eqref{eq:V} that the diagonal part
of~$\cV^o$ is zero.
\end{proof}

\begin{proof}[Proof of Lemma \ref{lem:UV}]\mbox{}\nopagebreak

$\eqref{enum:1}\Rightarrow\eqref{enum:2}$. If this system is
meromorphically isomorphic to the system with matrix $Udz/z^2$,
then the monodromy is clearly equal to identity. One also remarks
that the system is equivalent, by holomorphic base change, to a
system
\[
\Big(\frac Uz-V_{\diag}+zC(z)\Big)\frac{dz}z,
\]
where $V_{\diag}$ is the diagonal part of $V$ and $C(z)$ is
holomorphic. If it is holomorphically equivalent to $Udz/z^2$, let
us denote by $P_0+zP_1+\cdots$ a base change between both systems
(hence $P_0$ is invertible). Then $P_0$ is diagonal since it
commutes with~$U$, and therefore it commutes with $V_{\diag}$. We
must then have $V_{\diag}=[U,P_1]$, which implies $V_{\diag}=0$
since $[U,P_1]_{\diag}=0$.

\medskip
$\eqref{enum:2}\Rightarrow\eqref{enum:1}$. Assume that the
monodromy is equal to identity.

(a) After a base change by a matrix in $\mathrm{GL}_n(\CC[\![
z]\!])$, the system takes the normal form
\[\tag{$*$}
\Big(\frac Uz-V_{\diag}\Big)\frac{dz}z,
\]
hence it is of ``exponential type'', and can be described by
Stokes data consisting of two Stokes matrices $S_+,S_-$ (cf. for
instance \cite{B-J-L79}). More precisely, consider polar
coordinates $z=\rho e^{i\theta}$ and denote by $\cA$ the sheaf
on~$S^1$ consisting of germs $f$ of $C^\infty$ functions on
$\mathbb R_+\times S^1$ satisfying $\ov z\partial_{\ov z}f=0$ on
$\{\rho\neq0\}$. Then each base change in $\mathrm{GL}_n(\CC[\![
z]\!])$ as above can be lifted as a base change in
$\Gamma(I,\mathrm{GL}_n(\cA))$ so that the new matrix is $(*)$, if
$I$ is a small neighbourhood of a closed interval of $S^1$ of
length~$\pi$ with ``generic'' boundary points. The Stokes matrices
$S_+,S_-$ compute the multiplicative difference between base
changes corresponding to two opposite intervals at each of the
boundary points. In particular, if $S_+=\Id$ and $S_-=\Id$, then
there is a base change in $\Gamma(S^1,\mathrm{GL}_n(\cA))$. Since
$\Gamma(S^1,\cA)=\CC\{z\}$, the latter group is
$\mathrm{GL}(\CC\{z\})$. In other words, if $S_+=\Id$ and
$S_-=\Id$, then the system is holomorphically equivalent to a
system $(*)$.

(b) On the other hand, the monodromy of the system can be
presented as a product $S_+^{-1}S_-$, where, in a suitable basis,
the matrix of $S_+$ is upper triangular, and that of $S_-$ is
lower triangular. Therefore, if the monodromy is the identity, we
have $S_+=S_-$, so both are diagonal, and thus the formal
monodromy, which is equal to $S_{+,\diag}^{-1}S_{-,\diag}$, is
also equal to $\Id$, and the Stokes data is equivalent to the data
$S_+=\Id$, $S_-=\Id)$.

(c) Since the Stokes data are equal to identity, the system is
holomorphically equivalent to the associated formal system $(*)$.
Since the formal monodromy is the identity, $V_{\diag}$ has
integral entries. By a suitable rescaling of the basis by powers
of~$z$, the system is then meromorphically equivalent to a system
with $V_{\diag}=0$, that is, $Udz/z^2$.

(d) If $V_{\diag}$ is already zero, the rescaling is not necessary
and the isomorphism can be chosen holomorphic.

Assume now that \eqref{enum:1}, equivalently \eqref{enum:2}, is
satisfied. We already know from the previous proof that
$V_{\diag}$ has integral entries (resp.~is zero). We notice that,
in the neighbourhood of $\infty$ with coordinate $z'=1/z$, the
connection is written $(V-\nobreak z'U)dz'/z'$, hence has a simple
pole at $z'=0$. The last assertion of the lemma is then a
consequence of the following lemma.
\end{proof}

\begin{lemma}
Let $\nabla$ be a meromorphic connection with a simple pole at
$z'=0$. Let $(V+z'C(z'))dz'/z'$ be its matrix. Assume that the
monodromy is equal to identity. Then $V$ is semi-simple with
integral eigenvalues.
\end{lemma}

\begin{proof}
We will consider the Levelt normal form (see
e.g.~\cite[Ex.~II.2.20]{Sabb}). Denote by $D$ the diagonal matrix
whose entries are the integral parts of the real parts of the
eigenvalues of $V$, that we can assume to be ordered as
$\delta_1\geq\cdots\geq\delta_n$. We can assume that $V$ is
block-diagonal, the blocks corresponding to distinct eigenvalues
of~$V$. We first consider the block-indexing by the distinct
eigenvalues of $D$. Then (see loc.~cit.), the monodromy matrix can
be written as $\exp\big(-2\pi i(V-D+T)\big)$, where $T$ is
strictly block-lower triangular. That $\exp\big(-2\pi
i(V-D+T)\big)=\mathrm{Id}$ implies first that $T=0$. Now, $V-D$ is
block-diagonal with respect to distinct eigenvalues of $V$, and
each block $B$ satisfies $\exp(-2\pi i B)=\mathrm{Id}$, which
implies that~$B$ is semi-simple with integral eigenvalues. Hence
so is $V$.
\end{proof}

\begin{remark}\label{rem:nec}
As a consequence, a necessary condition to have a holomorphic
isomorphism $\phi$ (provided that $\cQ=0$) is that $\cV$ is
semi-simple with integral eigenvalues. With the notation of Lemma
\ref{lem:UV}, if $U$ is regular semi-simple and $V$ is semi-simple
with integral eigenvalues, there is an inductive procedure to
check whether the monodromy is the identity or not. But the
condition on $U,V$ is not easy to formulate. In order to compute
the monodromy, we consider the system at $z=\infty$, where it has
a regular singularity. Setting $z'=1/z$, the system has matrix
$(V+z'U)dz'/z'$. We then try to find a meromorphic base change
(locally with respect to the variable $z'$) so that the matrix of
the system is constant. Since the system has regular singularity,
such a base change is known to exist, and an inductive procedure
is known. Once the matrix is constant, it is easy to check whether
the monodromy is the identity or not. In rank two, it reduces to
the condition that the diagonal part of $V$ is zero: this is the
contents of the computation of Corollary \ref{Q02}.
\end{remark}

\begin{proof}[Proof of theorem \ref{thmA}]
In the chosen canonical coordinates, the matrix of the
endomorphism $\mathcal{U}$, denoted by $(U_{ij})$, is diagonal,
and
\begin{equation*}
(U_{ij})=\diag(u^1,\dots, u^m).
\end{equation*}

Let $\phi:\cH_1\to\cH_2$ be an isomorphism of holomorphic vector
bundles. Consider the local $C^\infty$ frame $S_i := (\pi^*
e_i)_{|_{\mathbb{C} \times M}}$ of $\mathcal{H}_1$ and
$\mathcal{H}_2$. It is $\ov\partial$-holomorphic, but not
$\ov\partial+z \Phi^{\dag}$-holomorphic, according to
\eqref{eq:Phidag} and Corollary \ref{dUdag}(1). Set
\begin{equation*}
\phi (S_i) = \sum_j \phi_{ij} \cdot S_j.
\end{equation*}
So $(\phi_{ij})$ is a non-degenerate $C^\infty$ matrix on
$\mathbb{C} \times M$.

\begin{claim}\label{claim:1}
A $C^\infty$ isomorphism $\phi:\cH_1\to\cH_2$ is holomorphic if
and only if the matrix
$(\psi_{ij}):=(\phi_{ij})\cdot\diag(\exp(-z\ov{\boldsymbol{u}}))$
is holomorphic.
\end{claim}

\begin{proof}
The component of $\ov\partial+z \Phi^{\dag}$ along $\partial_{\ov
z}$ is $d''_z$, so $\phi$ is holomorphic with respect to $z$ if
and only if $d''_z[\phi(S_i)]=0$ for all $i$. On the other hand,
$d''_z(S_j)=0$. Therefore this is equivalent to $\partial_{\ov
z}(\phi_{ij})=0$, or equivalently to $\partial_{\ov
z}(\psi_{ij})=0$.

On the other hand, the component of $\ov\partial+ z\Phi^\dag$ on
$\ov{e_j}$ is, according to Remark \ref{rem:phidag},
\begin{align*}
(\ov\partial+ z \Phi^\dag)_{\ov{e_j}} \phi(S_i)
& = \partial_{\ov u{}^j} \phi(S_i) + z \Phi^\dag_{\ov{e_j}} \phi(S_i) \\
& = \sum_k \partial_{\ov u{}^j} (\phi_{ik} S_k) +z \sum_k \phi_{ik} \Phi^{\dag}_{\ov{e_j}} S_k\\
& = \sum_k \partial_{\ov u{}^j} (\phi_{ik})\cdot S_k- z \phi_{ij}
S_j,
\end{align*}
and its vanishing is equivalent to
\begin{equation}\label{eq:holo2}
\begin{cases}
\partial_{\ov u{}^j} (\phi_{ik})=0,& \forall i,j, k \neq j,\\
\partial_{\ov u{}^j} (\phi_{ij})= z \phi_{ij},& \forall i, j,
\end{cases}
\end{equation}
that is, to $\partial_{\ov u{}^j} (\psi_{ik})=0$ for all $i,j,k$.
\end{proof}

\begin{claim}\label{claim:2}
The relation $\phi \circ\widetilde{\nabla}_{\partial_z} =
\widetilde{D}_{\partial_z} \circ\phi$ is equivalent~to
\begin{equation}\label{eq:CGF}
\partial_z (\psi)=\frac{1}{z^2} [\mathcal{U}, \psi] - \frac{1}{z} \cV \psi,
\end{equation}
with $\psi=\phi\cdot\exp(-z\cU^\dag)$.
\end{claim}

\begin{proof}
By straightforward computations, we get
\[
\tilde{\nabla}_{\partial_z} S_i=\frac{1}{z} \Big(\frac{1}{z}
\mathcal{U} - \cV\Big)(S_i)=\frac{1}{z} \Big[\frac{1}{z} u^i\cdot
S_i - \sum_k \cV_{ik} \cdot S_k\Big].
\]
So we have
\begin{align*}
\phi (\tilde{\nabla}_{\partial_z} S_i) & = \frac{1}{z}
\Big[\frac{1}{z} u^i\cdot \phi(S_i) - \sum_k
\cV_{ik} \cdot \phi(S_k) \Big] \\
& = \frac{1}{z^2} u^i\sum_j \phi_{ij} \cdot S_j - \frac{1}{z} \sum_{jk} \cV_{ik} \phi_{kj} \cdot S_j \\
& = \sum_j \Big[\frac{1}{z^2} u^i \phi_{ij} - \frac{1}{z} \sum_{k}
\cV_{ik} \phi_{kj} \Big] S_j.
\end{align*}
By similar computations, we get
\begin{align*}
\widetilde{D}_{\partial_z} \phi(S_i)
& = \widetilde{D}_{\partial_z} \Big(\sum_k \phi_{ik}\cdot S_k \Big) = \sum_k \partial_z (\phi_{ik})\cdot S_k + \sum_k \phi_{ik}\cdot \widetilde{D}_{\partial_z} S_k \\
& = \sum_k \partial_z (\phi_{ik})\cdot S_k + \frac{1}{z}\sum_k \phi_{ik}\cdot \Big[\frac{1}{z} \mathcal{U} - z \mathcal{U}^\dag \Big](S_k) \\
& = \sum_j \partial_z (\phi_{ij})\cdot S_j + \frac{1}{z}\sum_k \phi_{ik}\cdot \Big[\frac{1}{z} u^k S_k - z \ov{u^k}\cdot S_k \Big] \\
& = \sum_j \Big[ \partial_z (\phi_{ij})+ \frac{1}{z^2} \phi_{ij}
u^j -  \phi_{ij} \ov{u^j} \Big]S_j.
\end{align*}
The third equality holds because of Corollary \ref{dUdag}(1). So
the $\partial_z$ component of \eqref{eq:Cphi} is equivalent to
\begin{equation*}
\partial_z (\phi_{ij})+ \frac{1}{z^2}
\phi_{ij} u^j -\phi_{ij} \ov{u^j}= \frac{1}{z^2} u^i \phi_{ij} -
\frac{1}{z} \sum_{k} \cV_{ik} \phi_{kj}, \quad\forall i,j.
\end{equation*}
This amounts to the following, equivalent to \eqref{eq:CGF}:
\[
\partial_z (\phi_{ij}) - \phi_{ij}\ov u{}^j=\frac{1}{z^2} (u^i-u^j)\phi_{ij} - \frac{1}{z} \sum_{k} \cV_{ik} \phi_{kj}, \quad\forall i,j.
\qedhere
\]
\end{proof}

\subsubsection*{End of the proof of Theorem \ref{thmA}}
From Theorem \ref{Q0m}, the existence of $\phi$ is equivalent to
the existence of $\phi^o$ and, by the previous claims, to the
existence of a holomorphic invertible matrix $\psi^o$ satisfying
\ref{eq:thmA}. Notice also the entries of $\psi^o$ are entire
functions of $z$ which have moderate growth at infinity, since
\ref{eq:thmA} has a regular singularity at $z=\infty$. Therefore,
the entries of $\psi^o$ belong to $\CC[z]$. The same argument
applies to $(\psi^o)^{-1}$.
\end{proof}

\begin{proof}[Proof of Corollary \ref{CoFlka}]
By \cite{Lin09} (cf.~also Corollary \ref{dUdag}), we know that the
canonical coordinates $u^1,\dots, u^m$ are $\nabla$-flat (so that
$\nabla=D'$ is expressed as $d$ in the frame $(e_1,\dots,e_m)$)
and $\mathcal{Q}=0$. We also have $\cV_{ij}=0$ for all $i,j$, as a
consequence of \eqref{eq:V} and \cite[Eq.\,(63)]{Lin09}. According
to Claim 1 in the proof of Theorem \ref{thmA}, it is a matter of
showing that the holomorphic matrix $(\psi_{ij})$ is diagonal and
constant. The condition in Claim 2 above now reads
\[
\partial_z (\psi_{ij})=\dfrac{1}{z^2} (u^i-u^j)\psi_{ij}, \quad \forall i,j,\]
and the only holomorphic solutions consist of diagonal matrices
depending on $\boldsymbol u$ only. On the other hand, for every
$j$, the condition $\phi \circ \tilde{\nabla}_{e_j}
=\widetilde{D}_{e_j} \circ \phi$ is equivalent~to
\[
\partial_{u^j}(\phi_{il})=\frac{1}{z} (\delta_{jl}-\delta_{ij})\phi_{il}, \quad \forall i,l,
\]
or equivalently to
\[
\partial_{u^j}(\psi_{il})=\frac{1}{z} (\delta_{jl}-\delta_{ij})\psi_{il}, \quad \forall i,l,
\]
and, since $\psi$ is diagonal, it reduces to
$\partial_{u^j}(\psi_{ii})=0$, i.e., $\psi$ is constant.
\end{proof}

\begin{proof}[Proof of Corollary \ref{Q02}]
The condition $d/2\in\ZZ$ is necessary, as we have seen in Remark
\ref{rem:nec}. Let us show that it is sufficient. We will set
$n=d/2\in\ZZ$. We have $e=e_1+e_2$ and $g(e_1,e_2)=g(e_2,e_1)=0$.
Since we assume $g(e, e)=0$, we obtain
\[
\eta_2 := g(e_2, e_2)=g(e, e)-g(e_1, e_1)= -\eta_1.
\]
Recall also that $\eta_{12}=\partial_{u
^2}\partial_{u^1}\eta=\eta_{21}$. Therefore,
\[
\eta_{12}= \eta_{21} = -\eta_{11}= -\eta_{22},
\]
and from \eqref{eq:V} we obtain
\[
\cV_{12} = \frac{(u^2- u^1)\cdot \eta_{12}}{2 \eta_{2}}=
\frac{\mathcal{E}\eta_1}{2 \eta_{2}}= \frac{d}{2}=\cV_{21}.
\]
By Theorem \ref{thmA}, we are reduced to proving the existence of
$\psi^o\in\GL_2(\CC[z])$ such that
\[
z\partial_z\psi^o(z)=\frac1z\cdot[\cU^o,\psi^o(z)]-n\begin{pmatrix}0&1\\1&0\end{pmatrix}\psi^o(z).
\]
Let us set $\psi^o(z)=\sum_{k\geq0}\psi^o_kz^k$. The previous relation reduces to a recursive relation%
\refstepcounter{equation}\label{eq:rec}%
\begin{equation}\tag*{(\protect\ref{eq:rec})$_k$}\label{eq:reck}
[\cU^o,\psi^o_{k+1}]=\Big(n\begin{pmatrix}0&1\\1&0\end{pmatrix}+k\Id\Big)\psi^o_k
\end{equation}
with $\psi^o_k=0$ for $k<0$ and $k\gg0$, and $\psi^o_0$
invertible. Let us first notice that, if $\psi^o$ exists, then
$[\cU^o,\psi^o_0]=0$, that is, $\psi^o_0$ is diagonal. We will
show the existence and uniqueness of a solution $\psi^o$ with
$\psi^o_0=\Id$, and it will be clear that any solution will be of
the form $\psi^o\cdot\delta$, where $\delta$ is diagonal constant
and invertible.

Let us set
\[
A=\begin{pmatrix}0&1\\1&0\end{pmatrix},\quad
B=\begin{pmatrix}0&1\\-1&0\end{pmatrix},\quad
D=\begin{pmatrix}-1&0\\0&1\end{pmatrix}.
\]
Setting also $u^1_o-u^2_o=x$, we have
\[
[\cU^o,A]=xB,\quad [\cU^o,B]=xA,\quad AB=D,\quad AD=B.
\]
Then \ref{eq:reck} determines $\psi^o_{k+1}$ in terms of
$\psi^o_k$ up to a diagonal term $\delta_{k+1}$, which in turn is
determined by \eqref{eq:rec}$_{k+1}$ and the condition that
$[\cU^o,\psi^o_{k+1}]$ has zeros on the diagonal. One finds, for
$1\leq k\leq |n|$,
\[
\psi^o_k=\frac{\prod_{j=0}^{k-1}(j^2-n^2)}{k!\,x^k}\,\Big(\Id-\frac
kn A\Big)D^k.
\]
and $\psi^o_k=0$ for $k\geq |n|+1$.
\end{proof}

\providecommand{\bysame}{\leavevmode\hbox to3em{\hrulefill}\thinspace}
\providecommand{\MR}{\relax\ifhmode\unskip\space\fi MR }
\providecommand{\MRhref}[2]{%
  \href{http://www.ams.org/mathscinet-getitem?mr=#1}{#2}
}
\providecommand{\href}[2]{#2}


\begin{thebibliography}{10}

\bibitem{B-J-L79}
Werner Balser, Wolfgang~B. Jurkat, and Donald~A. Lutz, \emph{Birkhoff
  invariants and {S}tokes' multipliers for meromorphic linear differential
  equations}, J.~Math. Anal. Appl. \textbf{71} (1979), no.~1, 48--94.
  \MR{545861 (81f:34010)}

\bibitem{CV3}
Sergio Cecotti and Cumrun Vafa, \emph{Topological--anti-topological fusion},
  Nuclear Phys. B \textbf{367} (1991), no.~2, 359--461. \MR{MR1139739
  (93a:81168)}

\bibitem{CVN}
\bysame, \emph{On classification of {$N=2$} supersymmetric theories}, Comm.
  Math. Phys. \textbf{158} (1993), no.~3, 569--644. \MR{MR1255428 (95g:81198)}

\bibitem{D}
Boris Dubrovin, \emph{Geometry of {$2$}{D} topological field theories},
  Integrable systems and quantum groups (Montecatini Terme, 1993), Lecture
  Notes in Math., vol. 1620, Springer, Berlin, 1996, pp.~120--348.
  \MR{MR1397274 (97d:58038)}

\bibitem{Hert}
Claus Hertling, \emph{Frobenius manifolds and moduli spaces for singularities},
  Cambridge Tracts in Mathematics, vol. 151, Cambridge University Press,
  Cambridge, 2002. \MR{MR1924259 (2004a:32043)}

\bibitem{Hert2}
\bysame, \emph{{$tt\sp *$} geometry, {F}robenius manifolds, their connections,
  and the construction for singularities}, J. Reine Angew. Math. \textbf{555}
  (2003), 77--161. \MR{MR1956595 (2005f:32049)}

\bibitem{Lin09}
Jiezhu Lin, \emph{{Some constraints on Frobenius manifolds with a
  tt*-structure}}, Math. Z. \textbf{267} (2011), 81--108, online: 2009, DOI
  10.1007/s00209-009-0610-z.

\bibitem{Malg}
Bernard Malgrange, \emph{{Sur les d\'eformations isomonodromiques, I, II}},
  {S\'eminaire E.N.S. Math\'ematique et Physique} (L.~Boutet~{de Monvel},
  A.~Douady, and J.-L. Verdier, eds.), Progress in Math., vol.~37,
  Birkh{\"a}user, Basel, Boston, 1983, pp.~401--438.

\bibitem{Mani}
Yuri~I. Manin, \emph{Frobenius manifolds, quantum cohomology, and moduli
  spaces}, American Mathematical Society Colloquium Publications, vol.~47,
  American Mathematical Society, Providence, RI, 1999. \MR{MR1702284
  (2001g:53156)}

\bibitem{Sabb}
Claude Sabbah, \emph{{D\'eformations isomonodromiques et vari\'et\'es de
  Frobenius}}, Savoirs Actuels, CNRS~{\'E}ditions \& EDP~Sciences, Paris, 2002,
  English Transl.: Universitext, Springer \& EDP~Sciences, 2007.

\bibitem{Sabb08}
\bysame, \emph{{Universal unfoldings of Laurent polynomials and tt*
  structures}}, {From Hodge theory to integrability and TQFT: tt*-geometry}
  (R.~Donagi and K.~Wendland, eds.), Proc. Symposia in Pure Math., vol.~78,
  American Math. Society, Providence, RI, 2008, pp.~1--29.

\end{thebibliography}
\end{document}